\documentclass[11pt]{article}
\usepackage{graphicx}
\usepackage{tikz}
\usepackage{tikz-cd}
\usetikzlibrary{matrix}
\usetikzlibrary{intersections,shapes.arrows}
\usepackage{hyperref}
\usepackage{stmaryrd}
\usepackage{amssymb}
\usepackage{amsthm}
\usepackage{mathrsfs} 
\usepackage{amsmath}
\usepackage{enumitem, color, amssymb}
\usepackage[utf8]{inputenc} 
\usepackage[T1]{fontenc}    
\usepackage{url}     

\theoremstyle{definition}
\newtheorem{lem}{Lemma}[section]
\newtheorem{Prop}[lem]{Proposition}
\newtheorem{Theo}[lem]{Theorem}
\newtheorem{Rem}[lem]{Remark}
\newtheorem{Cor}[lem]{Corollary}

\newtheorem{Def}[lem]{Definition}

\newtheorem*{Int}{Theorem}

\newcommand\DistTo{\xrightarrow{
   \,\smash{\raisebox{-0.5ex}{\ensuremath{\scriptstyle\sim}}}\,}}
\newcommand{\bew}{\hfill $\boxempty$     
}
\newcommand{\Z}{\mathcal{Z}}
\newcommand{\R}{\mathrm{Reg}([0,1],\mathbb{R})}
\newcommand{\Ra}{\mathrm{Reg}^\circ([0,1],\mathbb{R})}
\newcommand{\Ro}{\mathrm{BReg}((0,1),\mathbb{R})}
\newcommand{\Rb}{\mathrm{Reg}((0,1),\mathbb{R})}
\newcommand{\Rbc}{\mathrm{Reg}^\circ((0,1),\mathbb{R})}
\newcommand{\Rc}{\mathrm{BReg}^\circ((0,1),\mathbb{R})}
\newcommand{\Ry}{\mathrm{Reg}([0,1),\mathbb{R})}
\newcommand{\Rz}{\mathrm{Reg}((0,1],\mathbb{R})}
\newcommand{\Rx}{\mathrm{Reg}_S([0,1],\mathbb{R})}

\title{Maximal Spectra of rings consisting of regulated functions}
\author{Philipp Jukic}
\date{}
\begin{document}
\maketitle
\begin{abstract}
In 1969 Höchster proved that for every quasi-compact T1-space $X$ we can find a commutative ring $R$ such that $X$ is homeomorphic to the maximal spectrum $\mathrm{Specm}(R)$ of $R$. This result implies the existence of a commutative ring $R$ that admits a non-Hausdorff and totally disconnected maximal spectrum $\mathrm{Specm}(R)$.
However, there has not been an explicit example of such a commutative ring yet.
We will construct a ring $A$ such that $\mathrm{Specm}(A)$ is not Hausdorff and has only the singleton sets and itself as connected components. 
\end{abstract}
\section{Introduction}
We are interested in the topological properties of the maximal spectrum $\mathrm{Specm}(R)$ of a commutative ring $R$, where $\mathrm{Specm}(R)$ is endowed with the Zariski-topology, i.e with the subspace induced from that of $\mathrm{Spec}(R)$. To be more precise, we want to study the relationship between maximal spectra that are either Hausdorff or totally disconnected.  
By \cite[Proposition 11]{z1} for each quasi-compact T1-space $X$ we can find a commutative ring $R$ such that $X$ and $\mathrm{Specm}(R)$ are homeomorphic. Thus, there are commutative rings $R$ such that $\mathrm{Specm}(R)$ is totally disconnected but not Hausdorff, since one can find a topological T1-space that has these properties. An example of such a topological space is $X=\{-\infty\}\cup \mathbb{N}\cup \{\infty\}$, where the topology is given by the topological sum of the two one-point compactifications $\mathbb{N}\cup \{\infty\}$ and $\mathbb{N}\cup \{-\infty\}$. For more information see \cite{z23}.\par\smallskip       
That just proved the existence of a commutative ring $R$ such that $\mathrm{Specm}(R)$ has the desired topological properties. Our aim is now to give an explicit commutative ring that is as close as possible in satisfying the the above properties.\par\smallskip
Before we give a short summary of the methods displayed in \cite[Proposition 11]{z1}, will first consider the ring of continuous functions $C([0,1],\mathbb{R})$ over $[0,1]$. Each maximal ideal of $\mathrm{Specm}(C^0([0,1],\mathbb{R}))$ corresponds to a point $x\in [0,1]$. 
We have a basis for the topology of $\mathrm{Specm}(C^0([0,1],\mathbb{R}))$ given by $\{D_M(f):f\in C^0([0,1],\mathbb{R})\}$, where $D_M(f)$ is the principal open set $D(f)$ intersected with $\mathrm{Specm}(C([0,1],\mathbb{R}))$. Each of the principal open sets $D_M(f)$ corresponds to an open set in $[0,1]$ and $D_M(f)$ is connected if and only if it corresponds to a subinterval in $[0,1]$. Thus, $\mathrm{Specm}(C^0([0,1],\mathbb{R}))$ is Hausdorff and not totally disconnected.\par\smallskip
Consider again the space $X=\{-\infty\}\cup \mathbb{N}\cup \{\infty\}$ and the ring $C^\ast(X)$ of bounded continuous functions $X\rightarrow \mathbb{R}$. By \cite[4.9, p. 57]{z5} the maximal spectrum $\mathrm{Specm}(C^\ast(\mathbb{N}\cup \{\infty\}))$ equals $\mathbb{N}\cup \{\infty\}$ as a set, i.e there is a one-to-one correspondence 
\begin{equation}
\mathbb{N}\cup \{\infty\}\rightarrow \mathrm{Specm}(C^\ast(\mathbb{N}\cup \{\infty\})),x\mapsto \mathfrak{m}_x,
\end{equation}
where $\mathfrak{m}_x=\{f\in C^\ast(\mathbb{N}\cup \{\infty\}):f(x)=0\}$. 
It is obvious that $\mathbb{N}\cup \{\infty\}$ and $\mathrm{Specm}(C^\ast(\mathbb{N}\cup \{\infty\}))$ are not only equal as sets but also as topological spaces. However, this is not true anymore when we consider $X$ and $\mathrm{Specm}(C^\ast(X))$. 
Let $f\in C^\ast(X)$ with $f(-\infty)=b$ and $f(\infty)=a$. For open neighborhoods $U_b,U_a$ of $a,b$ we get 
$$\begin{array}{c}
f^{-1}(U_{b})=X\backslash K_1\\
f^{-1}(U_{a})=X\backslash K_2,
\end{array}$$     
for some finite sets $K_1,K_2\subseteqq X$. The vanishing of $a$ implies the vanishing of $b$ and vice versa.  In other words, $\mathfrak{m}_\infty=\mathfrak{m}_{-\infty}$. This shows that $\mathrm{Specm}(C^\ast(X))$ is homeomorphic to $\mathrm{Specm}(C^\ast(\mathbb{N}\cup \{\infty\}))$ and not to $X$.\par\smallskip
Using the methods in \cite{z1} we can find the following workaround: Start with a quasi-compact non Hausdorff T1-space $X$ that is totally disconnected. Let $V$ denote the set of all continuous functions $X\rightarrow W$, where $W=\{0,1\}$ is endowed with the topology given by $\{\varnothing, \{0,1\}, \{0\}\}$. Furthermore, let $W^V$ denote the set of all mappings $V\rightarrow W$. Let $\mathrm{ev}:X\rightarrow W^V$ be the evaluation map. Take $X'$ to be the closure of $\mathrm{ev}(X)$ in $W^V$ with respect to the patch topology defined in \cite{z5}. Finally, we endow $X'$ with the relative product topology of $W^V$. By \cite[Theorem 8, p. 56]{z5} the space $X'$ is nothing more than the spectralification of $X$. Furthermore, every point of $X$ is closed in $X'$. On the other hand every closed point of $X'$ is contained in $X$. Theorem \cite[Thoerem 6, p. 51]{z5} shows that $X'$ is homeomorphic to the spectrum of the ring $R$ given by all continuous functions $X'\rightarrow \mathbb{R}$, where $\mathbb{R}$ is endowed with the discrete topology. Altogether, $X$ is homeomorphic to $\mathrm{Specm}(R)$.\par\smallskip
However, it is obvious that this approach will not allow us to describe the ring $R$ in an explicit manner: One problem, for example, is to determine the set $X'$. 
To produce rings that give us maximal spectra that are as 'close' as possible in being totally disconnected and not Hausdorff, we will turn our attention to the rings of regulated functions $\R$ on the interval $[0,1]$. 
The main result there is 
\begin{Int}
The maximal spectrum $\mathrm{Specm}(\R)$ is totally disconnected and Hausdorff.
\end{Int}
This result is also present in \cite{z22} but we will prove the same statement with a different approach.\par\smallskip
In section \ref{sx2}, we will therefore consider the maximal prime ideals of regulated functions on an open interval. In fact, this consideration will lead us to the construction of a ring that gives us a non-Hausdorff maximal spectrum that has just the singletons and itself as connected components. The results in this section can be seen as an extension of the results presented in \cite{z22}.\par\smallskip
In sections \ref{se5} and \ref{se6} we will consider the ring $\Ra$ that consists of all functions in $\R$ that have only a finite amount of discontinuities. While section \ref{se5} is about the introduction and basic properties of these new rings, section \ref{se6} deals with connections to ultrafilters. The main result of these two sections is the following statement: 
\begin{Int}
The space $\mathrm{Specm}(\Ra)$ is a profinite space.
\end{Int}
Finally, in section \ref{se7} we will show that all the maximal prime ideals of the ring $\R$ are real. 
For the ring $\Ra$ we even have that the real spectrum is homeomorphic to the prime spectrum: 
\begin{Int}
	The ring $\Ra$ is real closed. Furthermore, $\mathrm{Sper}(\Ra)$ and $\mathrm{Spec}(\Ra)$ are homeomorphic. At last, the set $\mathrm{Sperm}(\R)$ of closed points of $\mathrm{Sper}(\R)$ and $\mathrm{Specm}(\R)$ are equal as sets. 
\end{Int}
Furthermore, we give an example of a one dimensional ring that has a maximal spectrum with the same properties as the maximal spectrum of the ring constructed in section \ref{sx2}. The main result here is 
\begin{Int}
There is a semi-algebraic space $(X,\mathscr{O}_X)$ with the following properties:
	\begin{enumerate}[label=\alph*]
		\item The maximal spectrum of the global sections $\mathrm{Specm}(\mathscr{O}_X(X))$ has only itself and the singleton subsets as its connected components.
	\item  The maximal spectrum of the global sections $\mathrm{Specm}(\mathscr{O}_X(X))$ is not Hausdorff.
	\item The dimension of $\mathscr{O}_X(X)$ is $1$. 
	\item $\mathrm{Sperm}(\mathscr{O}_X(X))$ is homeomorphic to $\mathrm{Specm}(\mathscr{O}_X(X))$. 
	\end{enumerate}
\end{Int}
This result is entirely independent of the work done in the previous section. Hence it cannot be understood in the framework of \cite{z22} and therefore does not extend the results there.\par\smallskip
The rest of the introduction is dedicated to introduce some central definitions and statements.
\begin{Def}\label{C1.d1}
A ring $R$ is a called Gelfand if for any two distinct maximal right ideals $\mathfrak{m}$ and $\mathfrak{m}'$ we can find two elements $a\in R\backslash \mathfrak{m}$ and $a'\in R\backslash \mathfrak{m}'$ such that $aRa'=0$. 
\end{Def}
Things, however, get a lot simpler, since we are considering only commutative rings. In fact, we have the following simple lemma:
\begin{lem}\label{C1.l1}
If $R$ is a commutative Gelfand-ring, then the characterization in Definition \ref{C1.d1} is equivalent to the following characterization: For each pair of distinct maximal ideals $\mathfrak{m}$ and $\mathfrak{m}'$, there are ideals $I_1$ and $I_2$ of $R$ such that $I_1\nsubseteqq \mathfrak{m}$, $I_2\nsubseteqq \mathfrak{m}'$ and $I_1I_2=0$.
\end{lem}
\textbf{Proof}: Obvious.\bew
\begin{Prop}
Let $R$ be a commutative ring. Then the following statements are equivalent:
\begin{enumerate}[label=\alph*] 
\item $\mathrm{Spec}(R)$ is normal.
\item The ring $R$ is Gelfand.
\end{enumerate}
\end{Prop}
\textbf{Proof}: That is Lemma \ref{C1.l1} + \cite[p. 2]{z3}.\bew\\\\
More interestingly, the property of being a Gelfand-ring encodes some topological information about $\mathrm{Specm}(R)$.  
\begin{Prop}\label{C1.P1}
Let $R$ be a ring. Then the following statements are equivalent:
\begin{enumerate}[label=\alph*]
\item $R$ is a Gelfand-ring. 
\item If $V_1,\ldots,V_n\subseteqq \mathrm{Spec}(R)$ are closed with $V_1\cap\cdots\cap V_n=\varnothing$, then there are $c_1,\ldots,c_n\in R$ with $\overline{D(c_1)}\cap \ldots\cap \overline{D(c_n)}=\varnothing$. 
\item For all $a,b\in R$ with $\mathcal{V}(a)\cap \mathcal{V}(b)=\varnothing$ there are $c,d\in R$ with $\mathcal{V}(a)\subseteqq D(c)$, $\mathcal{V}(b)\subseteqq D(d)$ such that $D(c)\cap D(d)=\varnothing$.
\end{enumerate}
\end{Prop}
\textbf{Proof}: See \cite[Theorem 4.3, p.706]{z2}.\bew\\\\
Let $R$ be a semi-primitive commutative ring, i.e the Jacobson radical $J(R)$ equals $\{0\}$. Suppose that $R$ is Gelfand. Then we choose two distinct points $\mathfrak{m},\mathfrak{m}'\in\mathrm{Specm}(R)$. By statement (b) in Proposition \ref{C1.P1} we can find $c_1,c_2\in R$ such that $\mathcal{V}_M(\mathfrak{m})\subseteqq D_M(c_1)$, $\mathcal{V}_M(\mathfrak{m}')\subseteqq D_M(c_2)$ and $D_M(c_1)\cap D_M(c_2)=\varnothing$. Thus, $\mathrm{Specm}(R)$ is Hausdorff. Conversely, suppose that $\mathrm{Specm}(R)$ is Hausdorff. Then for any two distinct points $\mathfrak{m},\mathfrak{m}'$ in $\mathrm{Specm}(R)$ we can find ideals $I_1$ and $I_2$ of $R$ such that $\mathfrak{m}\in D_M(I_1)$, $\mathfrak{m}'\in D_M(I_2)$ and $D_M(I_1)\cap D_M(I_2)=\varnothing$. Since $R$ is semi-primitive this can only work if $I_1I_2=\{0\}$. Hence, $R$ is Gelfand. In fact, we proved a well known fact:
\begin{Prop}\label{C1.p2cx}
Let $R$ be a commutative ring. The maximal spectrum $\mathrm{Specm}(R)$ is Hausdorff if and only if $R/J(R)$ is a Gelfand-ring.
\end{Prop}  
The usefulness of Proposition \ref{C1.p2cx} is that all the information is included in the ring $R$ if $R$ is semi-primitive. Since the ring of continuous functions and the ring of piecewise-continuous functions is semi-primitive, we do not need to deal with quotients of rings. In the following all rings that appear will be commutative. The importance of Proposition \ref{C1.p2cx} manifests itself in the following statement that will serve us as one of our main tools for analyzing $\mathrm{Specm}(R)$:
\begin{Prop}\label{C2.p1}
	Let $R$ be a ring. The following conditions are equivalent:
	\begin{enumerate}[label=\alph*]
		\item $R$ is a Gelfand-ring and $\mathrm{Specm}(R)$ is totally disconnected.
		\item $R$ is clean, i.e every element of $R$ can be written as a sum of a unit and an idempotent element of $R$.
	\end{enumerate} 
\end{Prop}
\textbf{Proof}: See \cite[Theorem 1.1]{z4}.\bew
\section{Ring of regulated functions}\label{s2} 
In this section we will introduce the ring of regulated functions on the interval $[0,1]$ and write down the most important properties of its maximal spectrum.
\begin{Def}
Let $\R$ denote the set of all functions $f:[0,1]\rightarrow \mathbb{R}$ such that the following statements are always satisfied:
\begin{itemize}
\item For every $c\in [0,1)$ the limit $f(c+)=\lim_{x\rightarrow c,x>c}f(x)$ exists.
\item For every $c\in (0,1]$ the limit $f(c-)=\lim_{x\rightarrow c,x<c}f(x)$ exists.
\end{itemize}
The functions in $\R$ are called regulated functions. The set of regulated functions $\Rb$ on the open interval $(0,1)$ is defined in the same way. 
\end{Def}
\begin{Rem}
	Strictly speaking the ring $\Rb$ and its elements are not called regulated. One of the key properties of actual regulated functions is that they can be approximated by step functions with respect to the supremum norm. That is not the case with functions out of $\Rb$.  
\end{Rem}
\begin{Def}\label{def1}
	We define $P^+:\R\rightarrow \Ry$ to be the operator that sends $f\in \R$ to $P^+f\in \Ry$, where $P^+f(x)=f(x+)$ for every $x\in [0,1)$. In the same manner, $P^-:\R\rightarrow\Rz$ is defined by $P^-f(x)=f(x-)$ for $x\in (0,1]$. The operators $P^+:\Rb\rightarrow \Rb$ and $P^-:\Rb\rightarrow \Rb$ are defined in the same way.   
\end{Def}
We will use the notation that we established in Definition \ref{def1}  for certain subrings of $\R$ and $\Rb$. It will be always clear from the context which definition of $P^+$ resp. $P^-$ to consider. 
\begin{Def}
Let $x_1\in [0,1)$, $x_2\in (0,1]$ and $x_3\in [0,1]$. We define 
$$
\begin{array}{c}
\mathfrak{P}_{x_1}^+=\{f\in \R:f(x_1+)=0\}\\
\mathfrak{P}_{x_2}^-=\{f\in \R:f(x_2-)=0\}\\
\mathfrak{M}_{x_3}=\{f\in \R:f(x_3)=0\}
\end{array}
$$
If not otherwise stated, the same notation will also be used for other subrings of $\R$. 
\end{Def}    
\begin{Prop}\label{C2.lx}
Let $x_1\in [0,1)$, $x_2\in (0,1]$ and $x_3\in [0,1]$. The sets $\mathfrak{P}_{x_1}^+$, $\mathfrak{P}_{x_2}^-$ and $\mathfrak{M}_{x_3}$ are all maximal ideals in $\R$. Moreover, we have $$\mathrm{Specm}(\R)=\{\mathfrak{P}_{x_1}^+:x_1\in [0,1)\}\cup \{\mathfrak{P}_{x_2}^-:x_2\in (0,1]\}\cup \{\mathfrak{M}_{x_3}:x_3\in [0,1]\}.$$ 
\end{Prop}
\textbf{Proof}: See \cite[Theorem 5, p.20]{z22}.\bew 
\begin{lem}\label{lx1}
The idempotent elements in $\R$ are the characteristic functions $\chi_A$, where $A$ is a union of sub-intervals of $[0,1]$. 
\end{lem}
\textbf{Proof}: See \cite[Theorem 3, p.19]{z22}.\bew
\begin{Theo}\label{tx1}
The maximal spectrum $\mathrm{Specm}(\R)$ is totally disconnected and Hausdorff.
\end{Theo}
\textbf{Proof}:  According to Proposition \ref{C1.p2cx} it is enough to prove that the ring $\R$ is Gelfand. By Proposition \ref{C2.p1} it is enough to verify that $\R$ is clean. For every function $f\in \R$ we can find a characteristic function $\chi_A$ and some $\varepsilon>0$ such that for all $x\in [0,1]$ we have $|f(x)-\chi_A(x)|\geq \varepsilon$ for a suitable union $A$ of subintervals. By Lemma \ref{lx1} we know that $\chi_A$ is an idempotent an therefore we are done.\bew 
\begin{Rem}
	Proposition \ref{tx1} gives an alternative proof of \cite[Theorem 1, p.19]{z22} without using the fact that the step functions are dense in $\R$ with respect to the supremum norm. Furthermore, Theorem \ref{tx1} shows that $\mathrm{Specm}(\R)$ is a profinite space.    
\end{Rem}            
\section{Regulated functions on the open interval $(0,1)$}\label{sx2}
In the following we will consider the rings $\Rb$ and its three subrings:
\begin{itemize} 
\item $\Rbc$, the ring of all functions of $\Rb$ with finitely many discontinuities.  
\item $\Ro$, the ring of bounded regular functions on $(0,1)$.
\item $\Rc$, the ring of bounded regular functions on (0,1) that have only a finite number of discontinuities. 
\item In general, we will use the notation $\mathrm{Reg}^\circ$ to denote the subring that consists of all functions with finitely many discontinuities. 
\end{itemize}
Our primary concern will be the ring $\Rbc$. The major difference between $\R$ and $\Rbc$ is not that each function in $\Rbc$ has a finite number of discontinuities, but that all functions in $\Rbc$ are defined on the open interval $(0,1)$. 
Thus $\Rbc$ contains functions $f$ such that $f(0+)$ does not exist, which is impossible in the ring $\R$.\par\smallskip 
These new functions in $\Rbc$ give rise to a new set of maximal prime ideals.
Let $\mathfrak{Z}$ be a set of functions $f\in \Rbc$ that is maximal with respect to the following conditions:
\begin{enumerate}[label=\Alph*]
\item For any finite subset $L\subseteqq \mathfrak{Z}$ at least one of the sets $$\bigcap_{f\in L} \Z(f), \bigcap_{f\in L} \Z(P^+f),\bigcap_{f\in L} \Z(P^-f)$$ is not empty.  
\item We have     
$$\bigcap_{f\in \mathfrak{Z}} \Z(f)=\bigcap_{f\in \mathfrak{Z}} \Z(P^+f)=\bigcap_{f\in \mathfrak{Z}} \Z(P^-f)= \varnothing.$$ 
\end{enumerate}
\begin{Rem}
Let us take a closer look at the conditions (A) and (B):\par\smallskip
\underline{Condition A}: Condition A guarantees that the ideal generated by a finite amount of functions in $\mathfrak{Z}$ does not contain a unit.\par\smallskip
\underline{Condition B}: We need this condition to make sure that $\mathfrak{Z}$ is not contained in any maximal ideal that corresponds to the ones we already know about. More precisely, $\mathfrak{Z}$ is not contained in any of the ideals $\mathfrak{P}_x^+,\mathfrak{P}_x^-,\mathfrak{M}_x$ for $x\in (0,1)$.\par\smallskip
Conditions (A) and (B) impose obviously a restriction on the size of the sets $\mathfrak{Z}$. And this restriction is in such a way that $\mathfrak{Z}$ becomes a maximal ideal. 
\end{Rem}
\begin{lem}
Each set $\mathfrak{Z}$ maximal with respect to the conditions (A) and (B) is a maximal ideal in $\Rbc$, $\Rc$ and $\Rb$. Furthermore, $\mathfrak{Z}$ is not contained in any of the ideals $\mathfrak{M}_x,\mathfrak{P}_x^+$ and $\mathfrak{P}_x^-$ for $x\in (0,1)$.   
\end{lem}
\textbf{Proof}:
It is enough to prove the assertion for the ring $\Rbc$.\par\smallskip
\underline{Step 1: $\mathfrak{Z}$ is a proper ideal}: 
It is not difficult to see that $\mathfrak{Z}$ is closed under multiplication with $\Rbc$.   
Let $\langle \mathfrak{Z}\rangle$ be the ideal generated by $\mathfrak{Z}$. The set $\langle \mathfrak{Z}\rangle$ is a proper ideal since any sum $\sum_{i=1}^rf_i$ of elements $f_1,\ldots,f_r$ out of $\mathfrak{Z}$ is not a unit by condition (A). Furthermore, $\langle \mathfrak{Z}\rangle$ satisfies condition (A) and condition (B). By maximality of $\mathfrak{Z}$ we get $\mathfrak{Z}=\langle\mathfrak{Z}\rangle$.\par\smallskip
\underline{Step 2: The ideal $\mathfrak{Z}$ is not contained in $\mathfrak{P}_x^+,\mathfrak{P}_x^-$ or $\mathfrak{M}_x$}: That is condition (B).\par\smallskip
\underline{Step 3: The ideal $\mathfrak{Z}$ is maximal:} Any maximal ideal containing $\mathfrak{Z}$ must satisfy condition (A) and by step 2 also condition (B). The maximality of $\mathfrak{Z}$ with respect to the condition (A) and (B) implies that it is also a maximal ideal.\bew\par\smallskip
From condition (A) and (B) we get that $0\in \overline{\Z(f)}$ or $1\in \overline{\Z(f)}$ for any $f\in \mathfrak{Z}$. Thus the set $B$ of all the ideals $\mathfrak{Z}$ constructed above splits into a union $$B=B_1\cup B_0,$$ where $B_1$ consists of all $\mathfrak{Z}$ with $1\in \overline{\Z(f)}$ for all $f\in \mathfrak{Z}$ and $B_0$ of those with $0\in \overline{\Z(f)}$ for all $f\in \mathfrak{Z}$.  
The elements out of $B_1$ will be denoted with $\mathfrak{Z}_1$ and those out of $B_0$ will the denoted with $\mathfrak{Z}_0$.
\begin{Prop}\label{propid}
	The ideals of $\mathrm{Specm}(\Rbc)$ and $\mathrm{Specm}(\Rb)$ are given by the illustration below:
	\begin{center}
		\includegraphics[scale=0.5]{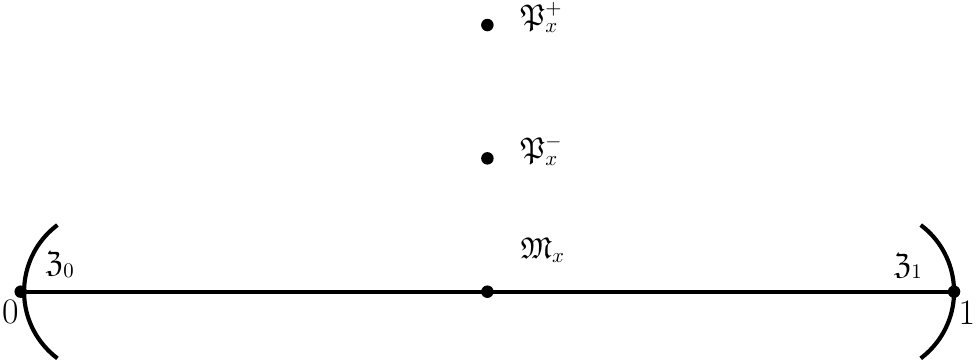}
	\end{center}
\end{Prop} 
\textbf{Proof}: Suppose that there is a maximal ideal $\mathfrak{P}$ that does not equal $\mathfrak{P}_x^+,\mathfrak{P}_x^-,\mathfrak{M}_x$ for all $x\in (0,1)$ . Then we know that for each $x\in (0,1)$ there are functions $f_1,f_2,f_3\in \mathfrak{P}$ such that the three inequalities  $$f_1(x)\neq 0, f_2(x+)\neq 0, f_3(x-)\neq 0$$ are satisfied. In other words, we have $$\bigcap_{f\in \mathfrak{P}} \Z(f)=\bigcap_{f\in \mathfrak{P}} \Z(P^+f)=\bigcap_{f\in \mathfrak{P}} \Z(P^-f)=\varnothing.$$ 
Note, that this resembles condition (B) from the construction of the ideals $\mathfrak{Z}$. Now, if $L\subseteqq \mathfrak{P}$ is some finite subset, then one of the intersections $\bigcap_{f\in L} \Z(f)$, $\bigcap_{f\in L} \Z(P^+f)$, $\bigcap_{f\in L} \Z(P^-f)$ is not empty. If all of these intersections would be empty, i.e $$\Z\left(\sum_{f\in L}f^2\right)= \Z\left(\sum_{f\in L}P^+f^2\right)= \Z\left(\sum_{f\in L}P^-f^2\right)=\varnothing ,$$ then $\sum_{f\in L}f^2$ is a unit.
Note, that this resembles condition (A). Now, it is not hard to see that $\mathfrak{P}$ is contained in $B_1$ or $B_0$.\bew  
\begin{Def}
	For every $n>>1$ we define 
	$\phi_n:(0,1)\rightarrow \mathbb{R}$ by $$\phi_n(x)=\begin{cases}
	1\,\text{if}\,x\geq \frac {1} {n}\\
	0\,\text{if}\,x<\frac {1} {n}\end{cases}.$$
	Let $I_0$ be the ideal generated by the family $(\phi_n)_{n\in \mathbb{N}}$. 
	Similarly, we define the ideal  
	$I_1$ to be generated by the family $(\phi_n')_{n\in \mathbb{N}}$, where 
	$$\phi_n'(x)=\begin{cases}
	1\,\text{if}\,x\leq 1-\frac {1} {n}\\
	0\,\text{if}\,x>1-\frac {1} {n}\end{cases}.$$
\end{Def}
\begin{Rem}
Obviously, $I_0$ and $I_1$ are not prime. 
\end{Rem}
\begin{Rem}\label{remul1}
	Let $V=\mathcal{V}_M(I)$ be a proper closed subset of the closed set $B_0$. Then there is a function $f\in I$ such that $\Z(f)\cap (0,\varepsilon)$ is an infinite discrete set for some small $\varepsilon>0$. That follows from the fact that each $f$ has only finitely many discontinuities and that $I$ must satisfy condition (A) and (B). To every such function we associate the set $\mathrm{Con}( \Z(f))$ of connected components of $\Z(f)$ and define the set 
	$$\mathcal{S}(I)=\{\mathrm{Con}( \Z(f)):f\in I,|\mathrm{Con}( \Z(f))|=\infty\}$$
	If $I$ is maximal, then it is easy to see that $\mathcal{S}(I)$ satisfies the following conditions:
	\begin{itemize}
		\item The empty set is not contained in $\mathcal{S}(I)$.
		\item The set $\mathcal{S}(I)$ is stable under intersections. 
		\item Any subset $S\subseteqq (0,1)$ that is connected and contains one element out of $\mathcal{S}(I)$ as a subset is already contained in $\mathcal{S}(I)$.  
	\end{itemize}
    Conversely, the three conditions give rise to a unique element in $B_0$. In other words, there is a 1-to-1 correspondence $\mathfrak{Z}_0\mapsto \mathcal{S}(\mathfrak{Z}_0)$. The same is also true if we consider $B_1$ instead of $B_0$. However, if we consider the ring $\Rb$ instead $\Rbc$ this correspondence does not hold anymore, since functions in $\Rb$ are allowed to have infinitely many jumps. In fact, this discussion shows that $\mathrm{Specm}(\Rb)$ and $\mathrm{Specm}(\Rbc)$ have prime ideals that are given by the same characterizations but the topological spaces $\mathrm{Specm}(\Rb)$ and $\mathrm{Specm}(\Rbc)$ fail to be homeomorphic.\par\smallskip
    The reason why we are looking at functions $f\in \Rbc$ that have an infinite amount of connected components emerges from the following observation: Consider the one of the maximal ideals $\mathfrak{Z}_0$ and suppose, for the sake of simplicity, that there is a function $f\in \mathfrak{Z}_0$ such that $|\mathrm{Con}( \Z(f))|=1$. Obviously, $\Z(f)$ cannot be just one point, since $\Z(P^+f)$ and $\Z(P^-f)$ are also finite and therefore $\mathfrak{Z}_0$ would not be a proper ideal anymore. Assume that $\Z(f)$ is one interval $C\subseteqq (0,1)$ with $0,1\notin\overline{C}$. Without loss of generality, we can assume that $C$ is a proper subset of $(0,1)$. Let $C'\subseteqq (0,1)$ be a closed interval such that $C\subseteqq C'$ and $\Z(P^+f),\Z(P^-f)\subseteqq C'$. Consider the restriction 
    $$\Rbc\rightarrow \chi_{C'}\Rbc$$
    and $\mathfrak{Z}_0$ under this map, i.e $\chi_{C'}\mathfrak{Z}_0$. Now, we have the following sequence 
    $$\chi_{C'}\Rbc\DistTo \mathrm{Reg}^\circ(C',\mathbb{R})\DistTo \Ra.$$
    Using this sequence to transport $\chi_{C'}\mathfrak{Z}_0$ into $\Rbc$ shows, by construction of $\mathfrak{Z}_0$, that $\chi_{C'}\mathfrak{Z}_0$ is not contained in any proper ideal. Thus, $\chi_{C'}\mathfrak{Z}_0$ must contain a unit under these conditions.  It is not difficult to deduce that $\mathfrak{Z}_0$ must contain a unit, contradicting the fact that $\mathfrak{Z}_0$ is assumed to be proper. 
    If the interval $C$ satisfies $0\in \overline{C}$, then situation is different because we cannot find a closed interval in $(0,1)$ that contains $C$. However, it is not difficult to see that any function $f$ with that property $0\in \overline{C}$ is contained in every $\mathfrak{Z}_0\in B_0$. This shows that the relevant functions that really characterize $\mathfrak{Z}_0$ are the ones where the vanishing set has an infinite amount of connected components. 
\end{Rem}
\begin{lem}\label{laemal1}
	We have $\mathcal{V}_M(I_0)=B_0$ and $\mathcal{V}_M(I_1)=B_1$ in $\mathrm{Specm}(\Rbc)$. Furthermore, $\bigcap_{\mathfrak{Z}_0\in B_0}\mathfrak{Z}_0=I_0$ resp. $\bigcap_{\mathfrak{Z}_1\in B_1}\mathfrak{Z}_1=I_1$.   
\end{lem}
\text{Proof}: 
It is enough to prove the assertion for $I_0$. 
Let us start with the inclusion $\mathcal{V}_M(I_0)\supseteqq B_0$: Since $\langle \mathfrak{Z}_0,I_0\rangle$ is a proper ideal for every $\mathfrak{Z}_0\in B_0$, we get that $\langle \mathfrak{Z}_0,I_0\rangle=\mathfrak{Z}_0$ resp. $I_0\subseteqq \mathfrak{Z}_0$. The other inclusion $\mathcal{V}_M(I_0)\subseteqq B_0$ follows from the fact that $\bigcap_{f\in I_0} \Z(f)=\bigcap_{f\in I_0} \Z(P^+f)=\bigcap_{f\in I_0} \Z(P^-f)=\varnothing$.\par\smallskip
By now, we know that $I_0\subseteqq \bigcap_{\mathfrak{Z}_0\in B_0}\mathfrak{Z}_0$. Let us verify the opposite inclusion. The ideal $I_0$ consists of all functions $f\in \Rbc$ that satisfy $(0,\varepsilon)\subseteqq \Z(f)$ for some $\varepsilon>0$. However, any function $f\in \mathfrak{Z}_0$ that does not satisfy $(0,\varepsilon)\subseteqq \Z(f)$ for some $\varepsilon>0$ must satisfy $\mathrm{Con}(\Z(f))\in \mathcal{S}(\mathfrak{Z}_0)$. But it is clear that we can find $\mathfrak{Z}_0'\in B_0$ with $\mathrm{Con}(\Z(f))\notin \mathcal{S}(\mathfrak{Z}_0')$. Thus, $\bigcap_{\mathfrak{Z}_0\in B_0}\mathfrak{Z}_0\subseteqq I_0$.\bew      
\begin{Prop}\label{Propw1}
	Any closed proper subset $V\subseteqq \mathrm{Specm}(\Rbc)$ of $B_0$ resp. $B_1$ with more than one element is not connected.  
\end{Prop}     
\textbf{Proof}: Without loss of generality we can assume that $V=\mathcal{V}_M(I)$ is infinite. Under these circumstances the set $\mathcal{S}(I)$ is not empty. Let $L\subseteqq (0,1)$ be such that $\mathrm{Con}(L)\in \mathcal{S}(I)$. Take two other sets $L_1,L_2\subseteqq (0,1)$ with $0\in \overline{L_1},\overline{L_2}$ and $\mathrm{Con}(L)=\mathrm{Con}(L_1)\uplus \mathrm{Con}(L_2)$. For every $\mathfrak{Z}_0\in V$ we have that either $\mathrm{Con}(L_1)\in \mathcal{S}(\mathfrak{Z}_0)$ or $\mathrm{Con}(L_2)\in \mathcal{S}(\mathfrak{Z}_0)$. Now, take two  distinct $\mathfrak{Z}_0,\mathfrak{Z}_0'\in V$ such that there are functions $f\in \mathfrak{Z}_0$ and $g\in \mathfrak{Z}_0'$ with $\mathrm{Con}(\Z(f))=\mathrm{Con}(L_1)$ and $\mathrm{Con}(\Z(g))=\mathrm{Con}(L_2)$. Consider the ideals $I_1=\langle I,f\rangle$ and $I_2=\langle I,g\rangle$. It is not hard to see that if $\mathrm{Con}(L_1)\in \mathcal{S}(\mathfrak{Z}_0)$ resp. $\mathrm{Con}(L_2)\in \mathcal{S}(\mathfrak{Z}_0)$, then $\mathfrak{Z}_0\in \mathcal{V}_M(I_1)$ resp. $\mathfrak{Z}_0\in \mathcal{V}_M(I_2)$. Altogether we get $V=\mathcal{V}_M(I_1)\uplus \mathcal{V}_M(I_2)$.\bew  
\begin{Prop}\label{theof1}
	The topological space $\mathrm{Specm}(\Rbc)$ has two connected components: The two sets $$\mathcal{V}_M(I_0)$$ and $$\mathcal{V}_M(I_1).$$
\end{Prop}
\textbf{Proof}: Both sets are closed and Proposition \ref{Propw1} implies that both are connected. Now, there are no other connected sets:
Without loss of generality, we can take a closed set $V=\mathcal{V}_M(I)$ such that $\mathcal{V}_M(I_0)\cup \mathcal{V}_M(I_1)\subsetneqq V$. Let $\phi_1,\phi_2:[0,1]\rightarrow \mathbb{R}$ be given by 
$$\phi_1(x)=\begin{cases}
0\,\text{if}\,x<\frac {1} {2}\\
1\,\text{if}\,x\geq \frac {1} {2}
\end{cases}$$ and
$$\phi_2(x)=\begin{cases}
0\,\text{if}\,x\geq \frac {1} {2}\\
1\,\text{if}\,x< \frac {1} {2}
\end{cases}.$$
Then $$V=\mathcal{V}_M(\langle I,\phi_1\rangle)\uplus \mathcal{V}_M(\langle I,\phi_2\rangle).$$ \bew  
\begin{Prop}\label{propequ1}
$\mathrm{Specm}(\Rbc)$ and $\mathrm{Specm}(\Rc)$ are not homeomorphic.
\end{Prop}
\textbf{Proof}: 
Let $M_\infty$ be the set of all functions $f$ in $\Rc$ such that $|f(x)|\rightarrow 0$ for $x\rightarrow 0$ or $x\rightarrow 1$. Let $\iota:\Rc\hookrightarrow \Rbc$ the canonical inclusion. We define $M_\infty'$ to be the set of all functions in $\iota(M_\infty)$ that satisfy $f\notin \mathfrak{P}$ for all $\mathfrak{P}\in \mathrm{Specm}(\Rbc)$.
Finally, we set $\tilde{M}_\infty=\iota^{-1}(M_\infty ')$ and introduce the ring $$R:=\Rc_{\tilde{M}_\infty}.$$
We will show that $\mathrm{Specm}(R)$ is homeomorphic to $\mathrm{Specm}(\Rc)$, which will automatically prove the assertion since $\mathrm{Specm}(R)$ and $\mathrm{Specm}(\Rbc)$ are not homeomorphic.\par\smallskip
Since $M_\infty'\subseteqq \Rbc^\times$ we get an embedding $R\hookrightarrow \Rbc$ which extends $\iota$. In other words, we get an induced mapping $$\Psi:\mathrm{Specm}(\Rbc)\rightarrow \mathrm{Specm}(R).$$
\underline{Step I: $\Psi$ is injective:} The only situations that needs checking is when we restrict $\Psi$ onto $B_0\cup B_1$. Without loss of generality, we can consider $\Psi$ on $B_0$. Take two distinct prime ideals $\mathfrak{Z}_0,\mathfrak{Z}_0'\in B_0$. By Remark \ref{remul1} we can find functions $f,g\in \Rbc$ such that $\mathrm{Con}(\Z(f))\in \mathcal{S}(\mathfrak{Z}_0)$, $\mathrm{Con}(\Z(f))\notin \mathcal{S}(\mathfrak{Z}_0')$, $\mathrm{Con}(\Z(g))\in \mathcal{S}(\mathfrak{Z}_0')$ and $\mathrm{Con}(\Z(g))\notin \mathcal{S}(\mathfrak{Z}_0)$. Furthermore, $f$ and $g$ can be chosen to be in $\Rc$. From there it is not difficult to see that the images of $\mathfrak{Z}_0$ and $\mathfrak{Z}_0'$ cannot be the same under $\Psi$.\par\smallskip
\underline{Step II: $\Psi$ is surjective}: It is enough to show that the ring extension $\Rbc|R$ is integral. It is enough to show that a function $f\in \Rbc$ with the two properties $|f(x)|\rightarrow \infty$ as $x\rightarrow 0$ and $0<|f(x)|<\infty$ as $x\rightarrow 1$ satisfies an integral equation. There is some $\varepsilon>0$ such that $\Z(f|_{(0,\varepsilon)})=\varnothing$. Let $g\in \Rc$ be given by 
$$g(x)=\begin{cases}
\frac {1} {f(x)}\,\text{if}\,x\in (0,\varepsilon)\\
1\,\text{otherwise}
\end{cases}.$$
It is not difficult to see that $g$ is contained in $\tilde{M}_\infty$. Since $$f^2-\frac {1} {g}\cdot f\in R,$$
we are done.\bew
\begin{Prop}\label{propbi1} The following statements hold:
\begin{enumerate}[label=\alph*] 
	\item We have the following bijections of sets: 
	$$\mathrm{Specm}(\R)\DistTo\mathrm{Specm}(\Ra)$$
	\item None of the maximal spectra of the rings
	$$\R,\Ra,\Rb,\Rbc,\Ro,\Rc$$ are homeomorphic to each other. 
\end{enumerate}	
\end{Prop}
\textbf{Proof}: 
(a): That is clear.\par\smallskip
(b):
\underline{Step I: The space $\mathrm{Specm}(\R)$:}
The only case that needs checking is if $\mathrm{Specm}(\R)$ is homeomorphic to $\mathrm{Specm}(\Ra)$.
Both $\R$ and $\Ra$ are clean rings. Let $u\in \R$ be a function such that $\Z(P^+u)$ is an infinite discrete subset of $[0,1]$. Then the set $\mathcal{V}_M(u)$ is infinite. But all closed subsets in $\mathrm{Specm}(\Ra)$ that are contained in $\{\mathfrak{P}_x^+:x\in [0,1)\}$ are finite. Thus, these two spaces cannot be homeomorphic.\par\smallskip
\underline{Step II: The space $\mathrm{Specm}(\Ra)$:} This case is clear by what we have done in the first step.\par\smallskip
\underline{Step III: The space $\mathrm{Specm}(\Rb)$:} By Remark \ref{remul1} we already know that the spaces $\mathrm{Specm}(\Rb)$ and $\mathrm{Specm}(\Rbc)$ are not homeomorphic. The only case that needs checking is the one if $\mathrm{Specm}(\Rb)$ is homeomorphic to $\mathrm{Specm}(\Ro)$. From the canonical inclusions we get the following diagram   
\begin{center}
	\begin{tikzpicture}
	\matrix (m) [matrix of math nodes,row sep=3em,column sep=4em,minimum width=2em]
	{
	\mathrm{Spec}(\Ro) & & \mathrm{Spec}(\Rb)\\ 
		& &\\
	\mathrm{Spec}(\Rc) & & \mathrm{Spec}(\Rbc) \\};
	\path[-stealth]
	(m-1-3) edge node [above,rotate=90,opacity=1] {} (m-1-1)
    (m-1-1) edge node [above,rotate=90,opacity=1] {} (m-3-1)
    (m-3-3) edge node [above,rotate=90,opacity=1] {} (m-3-1)
	(m-1-3) edge node [below,rotate=90,opacity=1] {} (m-3-3)
	;
	\end{tikzpicture}  
\end{center}
Replacing $	\mathrm{Spec}(\Ro)$ and $ \mathrm{Spec}(\Rb)$ with $\mathrm{Specm}(\Ro)$ and $ \mathrm{Specm}(\Rb)$ leads to the diagram
\begin{center}
	\begin{tikzpicture}
\matrix (m) [matrix of math nodes,row sep=3em,column sep=4em,minimum width=2em]
{
	\mathrm{Specm}(\Ro) & & \mathrm{Specm}(\Rb)\\ 
	& &\\
	\mathrm{Spec}(\Rc) & & \mathrm{Spec}(\Rbc) \\};
\path[-stealth]
(m-1-3) edge node [above,rotate=90,opacity=1]  [below,rotate=90,opacity=1] {\resizebox{2.5cm}{0.3cm}{$\approx$}} (m-1-1)
(m-1-1) edge node [above,rotate=90,opacity=1] {} (m-3-1)
(m-3-3) edge node [above,rotate=90,opacity=1] {} (m-3-1)
(m-1-3) edge node [below,rotate=90,opacity=1] {} (m-3-3)
;
\end{tikzpicture}  
\end{center}
It is not hard to deduce that the above diagram gives us a homeomorphism between $\mathrm{Specm}(\Rc)$ and $\mathrm{Specm}(\Rbc)$. But that contradicts Proposition \ref{propequ1}.\par\smallskip
\underline{Step IV: The space $\mathrm{Specm}(\Rbc)$:} By what we know already and Proposition \ref{propequ1} this case is clear.\par\smallskip
\underline{Step V: The spaces $\mathrm{Specm}(\Ro)$ and $\mathrm{Specm}(\Rc)$:} Both spaces are certainly not homeomorphic to each other. By what we already know, each one of these spaces cannot be homeomorphic to any other space mentioned in the formulation of this proposition.\bew 
\begin{Theo}\label{theox1}
The rings $\Rbc/I_0$ and $\Rbc/I_1$ have the following properties
\begin{enumerate}[label=\alph*]
\item Both rings are semi-primitive and not clean.
\item The spaces $\mathrm{Specm}(\Rbc/I_0)$ and $\mathrm{Specm}(\Rbc/I_1)$ are not Hausdorff. 
\item The topological spaces $\mathrm{Specm}(\Rbc/I_0)$ and $\mathrm{Specm}(\Rbc/I_1)$ have only themselves and the singleton subsets as their connected components. 
\end{enumerate}
\end{Theo}
\textbf{Proof}: 
 It is enough to prove the assertions for $\Rbc/I_0$.\par\smallskip
(a):
The spaces $$\mathrm{Spec}(\Rbc/I_0)$$ and $$\mathcal{V}(I_0)\subseteqq \mathrm{Spec}(\Rbc)$$
are homeomorphic.
By Lemma \ref{laemal1} we know that the intersection of all ideals in $\mathcal{V}_M(I_0)$ is contained in $I_0$. Therefore, the intersection of all ideals in $\mathrm{Specm}(\Rbc/I_0)$ equals $\{0\}$. Hence, $\Rbc/I_0$ is semi-primitive.\par\smallskip
Th ring $\Rbc/I_0$ is not clean: Take the residue class $\sin(x^{-1})+I_0$ of $\sin(x^{-1})$. There cannot be an idempotent $e+I_0$ such that $\sin(x^{-1})+e+I_0$ is a unit $u+I_0$: For every element $g\in u+I_0\subseteqq \Rbc$ we have either $u(0+)\in \mathbb{R}$ or $|u(0+)|=\infty$, which is not true for every element out of the right hand side.\par\smallskip
(b): Let $\mathfrak{Z}_0$ and $\mathfrak{Z}_0'$ be two distinct elements in $B_0$. Let $\mathfrak{Z}_0\in D_M(f)$ and $\mathfrak{Z}_0'\in D_M(g)$ be two principal open sets. The intersection $D_M(f)\cap D_M(g)$ consists of all primes that do not contain $f$ and $g$. But there are infinitely many of them: It easy to construct an ideal $I$ with $\mathcal{V}_M(I)\subseteqq B_0$ such that $\mathrm{Con}(\Z(f)),\mathrm{Con}(Z(g))\notin \mathcal{S}(I)$, for example by forcing $I$ to contain a function $h$ with $\mathrm{Con}(\Z(h))\cap \mathrm{Con}(\Z(f))=\mathrm{Con}(\Z(h))\cap \mathrm{Con}(\Z(g))=\varnothing$. Since $\mathcal{V}_M(I)$ is not empty, we see that $D_M(f)\cap D_M(g)$ is not empty. Thus, $\mathrm{Specm}(\Rbc/I_0)$ is not Hausdorff.\par\smallskip
(c): It is clear that $\mathrm{Specm}(\Rbc/I_0)$ is connected. By Proposition \ref{Propw1} the other connected components are the singleton sets of $\mathrm{Specm}(\Rbc/I_0)$.\bew      
\section{Finite amount of discontinuities}\label{se5}
\begin{Def}
Let $S\subseteqq [0,1]$ be a finite subset. Then we define $\Rx$ to be the set of all functions $f\in \R$ such that all points where $f$ is discontinuous is contained in $S$. 
\end{Def}
\begin{lem}
For each finite set $S\subseteqq [0,1]$ the set $\Rx$ is a ring. For a finite subset $A\subseteqq [0,1]$ we define the function $\psi_A:[0,1]\rightarrow \mathbb{R}$ with the property that $\psi_A(x)=1$ if and only if $\chi_A(x)=0$ and $\psi_A(x)=0$ if and only if $\chi_A(x)=1$ for $x\in [0,1]$. Whenever $A$ is a union of sub-intervals of $[0,1]$ we have $\psi_A\in \R$. If $A\subseteqq S$, then $\psi_A\in \Rx$.  
\end{lem}
\textbf{Proof}: We have that $\psi_A=\chi_{[0,1]\backslash A}$. Since $A$ is finite, the set $[0,1]\backslash A$ can be written a union of sub-intervals of $[0,1]$. Thus, $\psi_A\in \R$. The second assertion is obvious.\bew  
\begin{lem}\label{lx3}
Let $\iota^\ast:\mathrm{Spec}(\R)\rightarrow \mathrm{Spec}(\Rx)$ be induced by the inclusion $\iota:\Rx\hookrightarrow \R$. The restriction of $\iota^\ast$ onto $\mathrm{Specm}(\R)$ gives a bijection 
$$\mathrm{Specm}(\Rx)\DistTo \mathrm{Specm}(\R)\backslash \left(\{\mathfrak{P}_x^+,\mathfrak{P}_x^-:x\notin S\}\right).$$
\end{lem}
\textbf{Proof}: It is enough to verify that \begin{equation}\label{eq1}\mathrm{Specm}(\Rx)=\{\mathfrak{p}_x^+:x\in S\}\cup \{\mathfrak{p}_x^-:x\in S\}\cup \{\mathfrak{m}_x:x\in [0,1]\}\end{equation}
for 
$$
\begin{array}{c}
\mathfrak{p}_{x}^+=\{f\in \Rx:f(x+)=0\}\\
\mathfrak{p}_{x}^-=\{f\in \Rx:f(x-)=0\}\\
\mathfrak{m}_{x}=\{f\in \Rx:f(x)=0\}
\end{array}
$$
and $S\subseteqq (0,1)$. 
Obviously, $\mathrm{Specm}(\Rx)$ is at least as big as the right-hand side of Equation \ref{eq1}. Let $\mathfrak{m}$ be an arbitrary prime out of $\mathrm{Specm}(\mathrm{Reg}_S([0,1],\mathbb{R}))$ and $\mathfrak{M}$ the ideal generated by $\mathfrak{M}$ in $\R$. Suppose that $\mathcal{V}_M(\mathfrak{M})\neq \varnothing$. Let $\mathfrak{M}'\in \mathcal{V}_M(\mathfrak{M})$. 
Maximality implies $\mathfrak{m}=\mathfrak{M}\cap \mathrm{Reg}_S([0,1],\mathbb{R})=\mathfrak{M}'\cap \mathrm{Reg}_S([0,1],\mathbb{R})$
and therefore $\mathfrak{m}=\mathfrak{m}_x$ or $\mathfrak{m}=\mathfrak{p}_x^+$ or $\mathfrak{m}=\mathfrak{p}_x^-$. Suppose that $\mathcal{V}_M(\mathfrak{M})=\varnothing$, i.e the ideal $\mathfrak{M}$ that is generated by $\mathfrak{m}$ in $\R$ is not a proper ideal in $\R$. That means there are $f_1,\ldots,f_n\in \mathfrak{m}$ and $g_1,\ldots,g_n\in \R$ such that \begin{equation}\label{eq3}
\left|\sum_{i=1}^ng_i(x)f_i(x)\right|\geq \varepsilon
\end{equation}
for some fixed $\varepsilon>0$ and all $x\in [0,1]$. Now, there cannot be one $x\in [0,1]$ such that the functions $f_1,\ldots,f_r$ are contained in $\mathfrak{m}_x$. 
For any $\delta>0$ we can choose a point $x\in [0,1]$ such that    
\begin{equation}\label{eq2} \sum_{i=1}^nf_i^2(x)<\delta,\end{equation}
since $\mathfrak{m}$ is assumed to be proper. We already know that $\sum_{i=1}^nf_i^2$ does not have a root in $[0,1]$. Hence, there is only one possibility on how the above Inequality \ref{eq2} can be satisfied:
There is some $y\in S$ such that $f_1,\ldots,f_r\in \mathfrak{p}_y^+$ or $f_1,\ldots,f_r\in \mathfrak{p}_y^-$. 
But then inequality \ref{eq3} 
cannot hold, resulting in a contradiction. Hence, $\mathcal{V}_M(\mathfrak{M})$ cannot be empty, and we are done.\bew   
\begin{Rem}\label{remx1}
The space $\mathrm{Specm}(\Rx)$ is not totally disconnected: Let $S$ be a finite set with $x=\min(S)$ and take a function $f\in \Rx$ such that $\Z(f)$ is a proper subinterval of $[0,x]$. Then $\mathcal{V}_M(f)$ is a connected component of $\mathrm{Specm}(\Rx)$.   
\end{Rem}
\begin{Theo}
$\mathrm{Specm}(\Rx)$ is Hausdorff and not totally disconnected. 
\end{Theo}
\textbf{Proof}: The first assertion follows from Remark \ref{remx1}. The second assertion follows from the following facts: The singleton sets $\{\mathfrak{M}_x\},\{\mathfrak{P}_x^+\}$ and $\{\mathfrak{P}_x^-\}$ are closed-open for $x\in S$. 
There is a split exact sequence of $\mathbb{R}$-vector spaces
\begin{equation}\label{split1}
0\rightarrow \prod_S\mathbb{R}\rightarrow \mathrm{Reg}_S([0,1],\mathbb{R})\rightarrow \chi_{[0,1]\backslash S}\mathrm{Reg}_S([0,1],\mathbb{R})\rightarrow 0,
\end{equation}
where the first map is given by $$\prod_S\mathbb{R}\rightarrow \mathrm{Reg}_S([0,1],\mathbb{R}),(z_s)_{s\in S}\mapsto \sum_s\chi_{\{s\}}z_s$$ and the second one $$ \mathrm{Reg}_S([0,1],\mathbb{R})\rightarrow \chi_{[0,1]\backslash S)}\mathrm{Reg}_S([0,1],\mathbb{R})$$ by restricting the function onto $[0,1]\backslash S$. Let $S=\{x_1,\ldots,x_k\}$. Then  
\begin{equation}\label{split2}
\chi_{[0,1]\backslash S}\mathrm{Reg}_S([0,1],\mathbb{R})\cong \chi_{[0,x_1)}C^0([0,1])\times \prod_i\chi_{(x_i,x_{i+1})}C^0([0,1])\times \chi_{(x_k,1]}C^0([0,1]).
\end{equation}
From the isomorphism in \ref{split2} it is not hard to see that $\mathrm{Specm}(\chi_{[0,1]\backslash S}\mathrm{Reg}_S([0,1],\mathbb{R}))$ is Hausdorff and not totally disconnected. Now, using the sequence \ref{split1} and embed $\prod_S\mathbb{R}$ into $\mathrm{Reg}_S([0,1],\mathbb{R})$ we see that $\mathrm{Specm}(\mathrm{Reg}_S([0,1],\mathbb{R}))\cong \mathrm{Specm}(\chi_{[0,1]\backslash S}\mathrm{Reg}_S([0,1],\mathbb{R}))\uplus D_M(\prod_S\mathbb{R})\cong \mathcal{V}_M(\prod_S\mathbb{R})\uplus D_M(\prod_S\mathbb{R})$ is Hausdorff and not totally disconnected.\bew\par\smallskip
Let $\mathcal{G}\subseteqq 2^{[0,1]}$ be the set of all finite subsets of $[0,1]$. The set $\mathcal{G}$ is partially ordered by inclusion and for any two $S_1,S_2\in \mathcal{G}$ with $S_1\subseteqq S_2$ we have an inclusion $\mathrm{Reg}_{S_1}([0,1],\mathbb{R})\hookrightarrow \mathrm{Reg}_{S_2}([0,1],\mathbb{R})$. From this information we can construct the ring $\varinjlim_{S\in \mathcal{G}}\Rx$. 
Applying the $\mathrm{Spec}$-functor gives morphisms $$\mathrm{Spec}(\mathrm{Reg}_{S_2}([0,1],\mathbb{R}))\rightarrow \mathrm{Spec}(\mathrm{Reg}_{S_1}([0,1],\mathbb{R}))$$
and these morphisms map closed points onto closed points, leading to 
$$\mathrm{Specm}(\mathrm{Reg}_{S_2}([0,1],\mathbb{R}))\rightarrow \mathrm{Specm}(\mathrm{Reg}_{S_1}([0,1],\mathbb{R})).$$
From this information we can construct the space $\varprojlim_{S\in \mathcal{G}}\mathrm{Specm}(\Rx)$.
\begin{Theo}\label{apptheo}
	The topological space $$\varprojlim_{S\in \mathcal{G}}\mathrm{Specm}(\Rx)$$ is homeomorphic to $$\mathrm{Specm}(\Ra).$$ 
\end{Theo}   
\textbf{Proof}: By \cite[Lemma 4.6 ,p. 5]{z7} we get a homeomorphism $$\mathrm{Specm}\left(\varinjlim_{S\in \mathcal{G}}\Rx\right)\cong \varprojlim_{S\in \mathcal{G}}\mathrm{Specm}(\Rx).$$
Now, $\varinjlim_{S\in \mathcal{G}}\Rx\cong \Ra$ implies that $\varprojlim_{S\in \mathcal{G}}\mathrm{Specm}(\Rx)$ equals $\mathrm{Specm}(\Ra)$ as a set and as a topological space.\bew
 \section{Spectra of direct products of rings}\label{se6} 
Let $C^0([0,1])$ the set of all continuous functions $[0,1]\rightarrow \mathbb{R}$. In this situation \cite[Theorem 2.3, p. 25]{z5} provides us with the following interesting relationship:
\begin{Theo} 
\begin{enumerate}[label=\alph*]
\item If $I$ is an ideal in $C^0([0,1])$, then the family 
$$\Z(I)=\{\Z(f):f\in I\}$$ 
is a filter.
\item If $\mathcal{F}$ is a filter, then the family $\Z(\mathcal{F})=\{f:\Z(f)\in \mathcal{F}\}$ is an ideal in $C^0([0,1])$.
\end{enumerate}
\end{Theo}  
\textbf{Proof}: See \cite[Theorem 2.3, p. 25]{z5}.\bew
\begin{Cor}\label{Ch3.c1}
There is a 1-to-1 correspondence between the ultra-filters of $C^0([0,1])$ and $\mathrm{Specm}(C^0([0,1]))$.
\end{Cor}
\begin{Theo}\label{takak1}
Let $I$ be an infinite index set and $S$ the set of all functions $I\rightarrow \{F:F\,\text{a finite set of positive prime integers}\}$. In addition let $\Phi\in S$ denote the blank function $\Phi(i)=\varnothing$ for all $i\in I$. If $\mathcal{F}$ is an ultrafilter on some $\sigma\in S\backslash \{\Phi\}$, then the set $\Z(\mathcal{F})=\{a\in \prod_I\mathbb{Z}:\varrho_a\in \mathcal{F}\}$ is a maximal ideal of $\prod_I\mathbb{Z}$, where $\varrho_a$ is defined by $\varrho_a(i)=\{p\in \sigma(i):p\,\text{divides}\,a_i\}$. Furthermore, for every maximal prime ideal $\mathfrak{m}$ of $\prod_I\mathbb{Z}$, there is an ultrafilter $\mathcal{F}$ on some $\sigma \in S\backslash \{\Phi\}$ such that $\mathfrak{m}=\Z(\mathcal{F})$.
\end{Theo}
\textbf{Proof}: That is \cite[Theorem 3]{z6} and \cite[Theorem 2]{z6}.\bew\par\smallskip

\begin{Theo}\label{txx1}
Let $S$ be the set of all functions $[0,1]\rightarrow \mathbb{F}_2$ and let $\Phi$ denote the blank function as in Theorem \ref{takak1}. Let $\mathcal{F}$ be some ultrafilter on some $\sigma\in S\backslash \{\Phi\}$. Each such ultrafilter gives rise to maximal ideals 
$$\Z(\mathcal{F})=\left\{f\in \R:\varrho_f\in \mathcal{F}\right\},$$ 
$$\Z^+(\mathcal{F})=\left\{f\in \R:\varrho_f^+\in \mathcal{F}\right\},$$ and
$$\Z^-(\mathcal{F})=\left\{f\in \R:\varrho_f^-\in \mathcal{F}\right\},$$
where $\varrho_f$ corresponds to the set $\Z(f)\in S$, $\varrho_f^+$ to the set $\Z(P^+f)\in S$ and $\varrho_f^-$ to the set $\Z(P^-f)\in S$. Moreover, the maps $$\begin{array}{c}\mathcal{F}\mapsto \Z(\mathcal{F}),\\
\mathcal{F}\mapsto \Z^+(\mathcal{F}),\\
\mathcal{F}\mapsto \Z^-(\mathcal{F}),\end{array}$$
are bijections from the set of ultrafilters over some $\sigma\in S\backslash \{\Phi\}$ onto the set of maximal ideals $\{\mathfrak{M}_{x},x\in [0,1]\}$, $\{\mathfrak{P}_{x}^+,x\in [0,1)\}$, and $\{\mathfrak{P}_{x}^-,x\in (0,1]\}$.   
\end{Theo}  
\textbf{Proof}: 
Recall the definition of a ultrafilter $\mathcal{F}$ on some function $\sigma$:
\begin{enumerate}[label=\alph*]
\item We have $\Phi\notin \mathcal{F}$ and $\sigma\in \mathcal{F}$.
\item For two functions $\varrho_1,\varrho_2\in \mathcal{F}$ let $L_1$ and $L_2$ be the corresponding sets in $2^D$ of $\varrho_1$ resp. $\varrho_2$. Then $\varrho_1\wedge \varrho_2$ is the function that corresponds to the set $L_1\cap L_2$ and is contained in $\mathcal{F}$.
\item For two functions $\varrho_1,\varrho_2$ with $L_1\subseteqq L_2$ we write $\varrho_1\leq \varrho_2$. If $\varrho\leq \sigma$, then either $\varrho\in \mathcal{F}$ or $\sigma\backslash \varrho\in \mathcal{F}$, where $\sigma\backslash \varrho$ is the function we get from the complement of the corresponding sets.
\end{enumerate}
Obviously, the ideals given by $\Z(\mathcal{F})$ correspond to the ideals $\{\mathfrak{M}_x:x\in [0,1]\}$, 
the ideals given by $\Z^+(\mathcal{F})$ correspond to the ideals $\{\mathfrak{P}_x^+:x\in [0,1)\}$, and the ideals given by $\Z^-(\mathcal{F})$ correspond to the ideals $\{\mathfrak{P}_x^-:x\in (0,1]\}$. In other words, the mappings
$$\begin{array}{c}\mathcal{F}\mapsto \Z(\mathcal{F}),\\
\mathcal{F}\mapsto \Z^+(\mathcal{F}),\\
\mathcal{F}\mapsto \Z^-(\mathcal{F}),\end{array}$$
are all bijections.\bew
\par\smallskip 
Let $S=\{x_1,\ldots,x_n\}$ be a finite subset of $[0,1]$. Reconsider sequence \ref{split1} given by 
$$0\rightarrow \prod_S\mathbb{R}\rightarrow \mathrm{Reg}_S([0,1],\mathbb{R})\rightarrow \chi_{[0,1]\backslash S}\mathrm{Reg}_S([0,1],\mathbb{R})\rightarrow 0.$$
Applying the functor $\varinjlim_{S\in \mathcal{G}}$ leads again to a split exact sequence 
$$0\rightarrow \varinjlim_{S\in \mathcal{G}}\prod_S\mathbb{R}\rightarrow \varinjlim_{S\in \mathcal{G}}\mathrm{Reg}_S([0,1],\mathbb{R})\rightarrow \varinjlim_{S\in \mathcal{G}}\chi_{[0,1]\backslash S}\mathrm{Reg}_S([0,1],\mathbb{R})\rightarrow 0.$$
The homeomorphism $$\mathrm{Specm}\left(\varinjlim_{S\in \mathcal{G}}\chi_{[0,1]\backslash S}\mathrm{Reg}_S([0,1]\right)\cong \varprojlim_{S\in \mathcal{G}}\mathrm{Specm}\left(\chi_{[0,1]\backslash S}\mathrm{Reg}_S([0,1]\right)$$ 
shows that 
 $\mathrm{Specm}(\varinjlim_{S\in \mathcal{G}}\chi_{[0,1]\backslash S}\mathrm{Reg}_S([0,1])$ is profinite:\par\smallskip The space $\varprojlim_{S\in \mathcal{G}}\mathrm{Specm}\left(\chi_{[0,1]\backslash S}\mathrm{Reg}_S([0,1]\right)$ compact Hausdorff and it is not hard to see that it is even totally disconnected, i.e profinite. 
 The exact sequence above tells us that $\mathrm{Specm}(\varinjlim_{S\in \mathcal{G}}\mathrm{Reg}_S([0,1],\mathbb{R}))$ is also profinite. In fact, we have found an alternative proof of the statement:  
\begin{Theo}\label{propho1}
	The space $\mathrm{Specm}(\Ra)$ is profinite space.
\end{Theo}
\textbf{Proof}: The same arguments as in Theorem \ref{tx1}.\bew
\begin{Rem}
For any finite set $S\subseteqq [0,1]$ we have that $\mathrm{Spec}(\prod_S\mathbb{R})$ equals $\mathrm{Spec}(\prod_S\mathbb{K})$ for any field $\mathbb{K}$. This shows that we can replace the field $\mathbb{R}$ in the two identifications that we get from the two split exact sequences above. However, there is no chance that we can replace $\mathbb{R}$ by an arbitrary ring $R$. While any product of fields is a von Neumann regular ring, the product of arbitrary rings is not always von Neumann regular. Take $\prod_S\mathbb{Z}$ instead of $\prod_S\mathbb{K}$: The ring $\prod_S\mathbb{K}$ is von Neumann regular and has dimension $0$. But the dimension of $\prod_S\mathbb{Z}$ is certainly not $0$.\par\smallskip
The product $\prod_S\mathbb{Z}$ appears in the following context: 
Let $\mathrm{Reg}_S((0,1),\mathbb{Z})$ be the subring of $\Rb$ of $\mathbb{Z}$-valued functions that only have discontinuities at $S\subseteqq (0,1)$. For a discrete set $S$ the sequence \ref{split1} becomes a split exact sequence  
$$0\rightarrow \prod_S\mathbb{Z}\rightarrow \mathrm{Reg}_S((0,1),\mathbb{Z})\rightarrow \prod_S\mathbb{Z}\times \mathbb{Z}\rightarrow 0$$
of $\mathbb{Z}$-modules. By identifying $(0,1)$ with $\mathbb{R}$ we can consider the ring $\mathrm{Reg}_{\mathbb{N}}(\mathbb{R},\mathbb{Z})$. 
In order to study $\mathrm{Specm}(\mathrm{Reg}_{\mathbb{N}}(\mathbb{R},\mathbb{Z}))$ we need to study the maximal ideals of a infinite copy of $\mathbb{Z}$. 
 Theorem \ref{takak1} tells us that the maximal ideals of $\mathrm{Reg}_{\mathbb{N}}(\mathbb{R},\mathbb{Z})$ correspond to certain ultrafilters. 
Theorems \ref{txx1} and \ref{takak1} illustrate how far more complicated $\mathrm{Specm}(\mathrm{Reg}_{\mathbb{N}}(\mathbb{R},\mathbb{Z}))$ is than $\mathrm{Specm}(\Rbc)$ resp. $	\mathrm{Specm}(\Rb)$.       
\end{Rem} 
\section{From a real algebraic point of view}\label{se7}
Finally, we are going to investigate the rings $\R$, $\Rbc$ and $\Ra$ about their real algebraic properties: In particular we are interested in the real prime ideals of these rings and in the real spectrum $\mathrm{Sper}$ of these rings. For more information about real prime ideals and the real spectrum we refer to \cite{z19}. 
\begin{Prop}\label{plas1}
	Every ideal in $\mathrm{Specm}(\R)$ is real.
\end{Prop}
\textbf{Proof}: Take $f_1,\ldots,f_r\in \R$ such that $f_1^2+\cdots +f_r^2\in \mathfrak{M}_x$. From $$\Z\left(f_1^2+\cdots +f_r^2\right)=\bigcap_{i=1}^r\Z(f_i)\supseteqq \{x\}$$ we get that $f_1,\ldots,f_r\in \mathfrak{M}_x$. To finish the proof it suffices to consider the case where the functions $f_1,\ldots,f_r\in \R$ satisfy $f_1^2+\cdots +f_r^2\in \mathfrak{P}_x^+$. By applying the operator $P^+$ we get 
$$\Z\left(P^+f_1^2+\cdots +P^+f_r^2\right)=\bigcap_{i=1}^r \Z(P^+f_i)\supseteqq \{x\},$$
implying that $f_1,\ldots,f_r\in \mathfrak{P}_x^+$.\bew
\begin{Prop}\label{plas2}
	Every ideal in $\mathrm{Specm}(\Rbc)$ is real.
\end{Prop}
\textbf{Proof}: The ideals $\mathfrak{P}_x^+,\mathfrak{P}_x^-,\mathfrak{M}_x\in \mathrm{Specm}(\Rbc)$ can be treated with the same arguments as in Proposition \ref{plas1}. It remains to verify that $\mathfrak{Z}_0$ is real. If $f_1^2+\cdots +f_r^2\in \mathfrak{Z}_0$, then $\Z(f_1^2+\cdots +f_r^2)=\bigcap_{i=1}^r \Z(f_i)\neq \varnothing$. Since  $\langle\mathfrak{Z}_0,f_1,\ldots,f_r\rangle$ is a again a proper ideal, we get that it equals $\mathfrak{Z}_0$ by maximality.\bew
\begin{Prop}
	For every $\mathfrak{P}\in \mathrm{Specm}(\R)$ the field $\R/\mathfrak{P}$ has the following properties:
	\begin{enumerate}[label=\alph*]
	\item If $\mathfrak{P}=\mathfrak{M}_x$ for some $x\in [0,1]$, then $\R/\mathfrak{P}\cong \mathbb{R}$.
	\item If $\mathfrak{P}=\mathfrak{P}^+_x$ for some $x\in [0,1)$, then $\R/\mathfrak{P}\cong \mathbb{R}$.
	\item If $\mathfrak{P}=\mathfrak{P}^-_x$ for some $x\in (0,1]$, then $\R/\mathfrak{P}\cong \mathbb{R}$. 
	    \end{enumerate}
\end{Prop}
\textbf{Proof}: (a): That is obvious.\par\smallskip
(b): Let $f\in \R$ be a function such that $f$ is continuous at $x$. Then it is obvious that $f$ can be written as a sum of an element in $\mathfrak{P}_x^+$ and one number in $\mathbb{R}$. Suppose that $f$ has a discontinuity at $x$. If we take $a=\lim_{c\rightarrow x,c>x}f(c)$, then $f-a\in \mathfrak{P}_x^+$, implying that $f$ is a sum of an element in $\mathfrak{P}_x^+$ and some number in $\mathbb{R}$. Altogether we have $\R/\mathfrak{P}_x^+\cong \mathbb{R}$.\par\smallskip
(c): The same argument as in (b).\bew  
\begin{Def}\label{defreal}
	A ring $R$ is called real closed if the following conditions are satisfied:
	\begin{enumerate}[label=\alph*]
	\item The set of squares of $R$ is the set of non-negative elements of a partial order $\leq$ on $R$. Furthermore, the ring $R$ together with the partial order $\leq$ is an f-ring.
	\item Let $a,b\in R$. If $0\leq a\leq b$, then $b|a^2$.
	\item For every $\mathfrak{p}\in \mathrm{Spec}(R)$ the ring $R/\mathfrak{p}$ is integrally closed and the field of fractions of $R/\mathfrak{p}$ is real closed. 	
	\end{enumerate}
\end{Def}
\begin{lem}\label{tf1}
	Let $S\subseteqq [0,1]$ be a finite subset. Then $\mathrm{Reg}_S([0,1],\mathbb{R})$ is a real closed ring.
\end{lem}
\textbf{Proof}: We will prove the assertion in the following two steps:\par\smallskip
\underline{Step I: The direct product $R_1\times R_2$ of two real closed rings is real closed:} Conditions (a) and (b) from Definition \ref{defreal} are clear. Only condition (c) needs checking. For each $\mathfrak{p}\in \mathrm{Spec}(R_1\times R_2)$ there is an $\mathfrak{p}_1\in \mathrm{Spec}(R_1)$ or $\mathfrak{p}_2\in \mathrm{Spec}(R_2)$ such that we have one of the two isomorphism $$\begin{array}{c}
(R_1\times R_2)/\mathfrak{p}\cong R_1/\mathfrak{p}_1\\
(R_1\times R_2)/\mathfrak{p}\cong R_2/\mathfrak{p}_2.\end{array}$$
From this it is clear that also condition (c) holds, i.e $R_1\times R_2$ is a real closed ring.\par\smallskip
\underline{Step II: The ring $\chi_{[0,1]\backslash S}\mathrm{Reg}_S([0,1],\mathbb{R})$ is real closed:} Without loss of generality, we can assume that $S=\{x_1,\ldots,x_n\}\subseteqq (0,1)$. We need to verify that the rings of the form $$\chi_{[0,x_1)}C^0([0,1]),\chi_{(x_1,x_2)}C^0([0,1]),\ldots,\chi_{(x_{n-1},x_n)}C^0([0,1]),\chi_{(x_n,1]}C^0([0,1])$$
are real closed.
We get a diagram of isomorphism 
\begin{center}
	\begin{tikzpicture}
	\matrix (m) [matrix of math nodes,row sep=3em,column sep=4em,minimum width=2em]
	{
		\chi_{(x_1,x_2)}C^0([0,1]) & &\dots & & \chi_{(x_{n-1},x_n)}C^0([0,1])\\ 
		& & \chi_{(0,1)}C^0([0,1])& &\\
	\chi_{[0,x_1)}C^0([0,1]) & &C^0([0,1]) & & \chi_{(x_n,1]}C^0([0,1]) \\};
	\path[-stealth]
	(m-1-1)  edge node [above,rotate=90,opacity=1]  [below,rotate=90,opacity=1] {\resizebox{2.5cm}{0.3cm}{$\approx$}} (m-1-3)
	(m-1-1)  edge node [above,rotate=90,opacity=1]  [above,rotate=70,opacity=1] {\resizebox{2.5cm}{0.3cm}{$\approx$}} (m-2-3)
	(m-1-5)  edge node [above,rotate=90,opacity=1]  [above,rotate=110,opacity=1] {\resizebox{2.5cm}{0.3cm}{$\approx$}} (m-2-3)
	(m-1-3)  edge node [above,rotate=90,opacity=1]  [below,rotate=90,opacity=1] {\resizebox{2.5cm}{0.3cm}{$\approx$}} (m-1-5)
	(m-3-1)  edge node [above,rotate=90,opacity=1]  [below,rotate=90,opacity=1] {\resizebox{2.5cm}{0.3cm}{$\approx$}} (m-3-3)
	(m-3-5)  edge node [above,rotate=90,opacity=1]  [below,rotate=90,opacity=1] {\resizebox{2.5cm}{0.3cm}{$\approx$}} (m-3-3)
    (m-2-3)  edge node [above,rotate=90,opacity=1]  [below,rotate=0,opacity=1] {\resizebox{1.5cm}{0.3cm}{$\approx$}} (m-3-3);
	\end{tikzpicture}  
\end{center}
which proves that these rings are all real closed.\par\smallskip
\underline{Step III: The ring $\mathrm{Reg}_S([0,1],\mathbb{R})$ is real closed:} By the split exact sequence \ref{split1} we know that $$\mathrm{Reg}_S([0,1],\mathbb{R})\cong \prod_S\mathbb{R} \times\chi_{[0,1]\backslash S}\mathrm{Reg}_S([0,1],\mathbb{R})$$
as a $\mathbb{R}$-vector space. But it is easy to verify that the mapping 
$$\prod_S\mathbb{R} \times\chi_{[0,1]\backslash S}\mathrm{Reg}_S([0,1],\mathbb{R})\rightarrow \mathrm{Reg}_S([0,1],\mathbb{R}),((z_s)_{s\in S},f)\mapsto \chi_{[0,1]\backslash S}f+\sum_s\chi_{\{s\}}z_s$$  
is an isomorphism of $\mathbb{R}$-algebras. Step II implies that $\mathrm{Reg}_S([0,1],\mathbb{R})$ is real closed.\bew 
\begin{Theo}\label{tls1}
	The ring $\Ra$ is real closed. Furthermore $\mathrm{Sper}(\Ra)$ and $\mathrm{Spec}(\Ra)$ are homeomorphic. At last, the set $\mathrm{Sperm}(\R)$ of closed points of $\mathrm{Sper}(\R)$ and $\mathrm{Specm}(\R)$ are equal as sets.  
\end{Theo}
\textbf{Proof}: From \cite[p. 34]{z30} we know that a direct limit of real closed rings is again real closed. This fact combined with Lemma \ref{tf1} proves the first assertion. The second assertion follows from the well known fact that the real spectrum of a real closed ring is homeomorphic to its prime spectrum. Finally, the rest follows from Proposition \ref{propbi1}.\bew   
\begin{Cor}
	The space $\mathrm{Sperm}(\R)$ is profinite.
\end{Cor}
\textbf{Proof}: That is Theorem \ref{tls1} and Theorem \ref{propho1}.\bew 
\begin{Def}
(See \cite[Definition 1, p. 1]{z25})  Let $\mathrm{Sper}(R)$ be the real spectrum of a ring $R$. A pair $(X,\mathscr{O}_X)$ consisting of
a constructible subset $X$ of $\mathrm{Sper}(R)$ and the sheaf $\mathscr{O}_X$ of abstract semialgebraic functions on $X$ is called a semialgebraic subspace of $\mathrm{Sper}(R)$. 
\end{Def}
\begin{Theo}
	There is a semi-algebraic space $(X,\mathscr{O}_X)$ with the following properties:
	\begin{enumerate}[label=\alph*]
		\item The maximal spectrum of the global sections $\mathrm{Specm}(\mathscr{O}_X(X))$ has only itself and the singleton subsets as its connected components.
		\item  The maximal spectrum of the global sections $\mathrm{Specm}(\mathscr{O}_X(X))$ is not Hausdorff.
		\item The dimension of $\mathscr{O}_X(X)$ is $1$. 
		\item $\mathrm{Sperm}(\mathscr{O}_X(X))$ is homeomorphic to $\mathrm{Specm}(\mathscr{O}_X(X))$. 
	\end{enumerate}
\end{Theo}
\textbf{Proof}: Let $X=\mathrm{Sper}(\mathbb{R}[\mathrm{x}])$. We will show that $(X,\mathscr{O}_X)$ satisfies all three conditions.\par\smallskip
(a): By \cite[Proposition 7.3.2, p. 146]{z19} we have that $\mathscr{O}_X(X)$ is the ring of all semi-algebraic functions $\mathbb{R}\rightarrow \mathbb{R}$. Let $f\in \mathscr{O}_X(X)$ be a semi-algebraic function that vanishes nowhere. The semi-algebraic mappings 
$\mathbb{R}\rightarrow \mathbb{R}^2,x\mapsto (x,f(x))$, $\mathbb{R}^2\rightarrow \mathbb{R},(x,y)\mapsto y$ and $\mathbb{R}\backslash \{0\}\rightarrow \mathbb{R}^2,x\mapsto (x,x^{-1})$ can be composed to a sequence $$\mathbb{R}\rightarrow \mathbb{R}^2\rightarrow \mathbb{R}\backslash \{0\}\rightarrow \mathbb{R}^2,$$ which gives a semi-algebraic mapping $\mathbb{R}\rightarrow \mathbb{R}^2,x\mapsto (x,x)$. Thus, every non-vanishing function in $\mathscr{O}_X(X)$ is a unit. It is easy to see that all the maximal prime ideals of $\mathscr{O}_X(X)$ are ideals of the form $\mathrm{m}_x=\{f\in \mathscr{O}_X(X):f(x)=0\}$. Since any semi-algebraic function in $\mathscr{O}_X(X)$ has only finitely many roots, we see that the non-trivial closed sets of $\mathrm{Specm}(\mathscr{O}_X(X))$ are all finite, which implies the assertion.\par\smallskip
(b): Follows from the fact that all the non-trivial closed sets of $\mathrm{Specm}(\mathscr{O}_X(X))$ are all finite.\par\smallskip
(c): That is \cite[Theorem 1.1, p.754]{z29}.\par\smallskip
(d): That is clear.\bew

\end{document}